\documentclass[11pt]{article}                         

\makeatletter
\def\draft{\textheight=10.5truein \textwidth=7.5truein \parindent=8pt
           \voffset=-1truein \topmargin=0Truein
           \ifcase \@ptsize \hoffset=-1.5truein \or \hoffset=-1.35truein
                        \or \hoffset=-1.15truein \fi}
\def\quality{\textheight=240mm \textwidth=160mm \topmargin=0Truein
             \ifcase \@ptsize \hoffset=-23mm \or \hoffset=-20mm
                          \or \hoffset=-15mm \fi}
\def\USdraft{\textheight=9.85truein \textwidth=7.5truein \parindent=8pt
             \voffset=-1.0truein \topmargin=0Truein
             \ifcase \@ptsize \hoffset=-1.5truein \or \hoffset=-1.35truein
                          \or \hoffset=-1.15truein \fi}
\def\USquality{\textheight=250mm \textwidth=170mm \topmargin=0Truein
               \voffset=-1truein
               \ifcase \@ptsize \hoffset=-23mm
                       \or \hoffset=-20mm \or \hoffset=-15mm \fi}
\makeatother 
\USquality

   \newcommand\?[1]{} 

\def\n{\noindent}  
\def\proof{\smallskip \noindent {\bf Proof. \ }}       
\def\blanksquare{\,\,\,$\sqcup\!\!\!\!\sqcap$}         
\def\qed{\hfill\blanksquare\linebreak\smallskip\par}   

\def\bline(#1,#2)(#3,#4)(#5){\put(#1,#2){\line(#3,#4){#5}}}  

\def\thname{Theorem}     \def\lmname{Lemma}      \def\prname{Proposition}
\def\dfname{Definition}  \def\crname{Corollary}  \def\rmname{Remark}

\newtheorem{theorem}{\thname}[section]   
\newtheorem{lemma}[theorem]{\lmname}     

\newtheorem{dftn}{\dfname}[section]
\newtheorem{rmrk}[theorem]{\rmname}
\newenvironment{remark}{\begin{rmrk}\rm}{\end{rmrk}}     

             \catcode`@=11 \@addtoreset{equation}{section} \catcode`@=12
\newcommand\mlbscale{1pt} 
\newif\iffigs\figstrue 
\def\Bfig(#1,#2)#3#4{\begin{figure} \begin{center}
    \setlength{\unitlength}{\mlbscale}
       \iffigs \begin{picture}(#1,#2) #3 \end{picture}
       \else \begin{picture}(60,10)(0,0)
                   \put(0,0){\framebox(60,10){Figure}} \end{picture} \fi
    \end{center} \caption{#4} \end{figure}}

\def\bpic(#1,#2)#3{\setlength{\unitlength}{\mlbscale}
    \begin{picture}(#1,#2) #3 \end{picture}}
\def\IR{\hbox{\rm I\kern-.2em\hbox{\rm R}}} \def\ep{\varepsilon}
\def\IZ{\hbox{{\rm Z}\kern-.3em{\rm Z}}}

\def\map{T}

   \def\la{\lambda}

\def\noprint#1{}

  \def\map{T}  \def\noprint#1{}

\def\cX{{\cal X}}

\def\blim#1#2{\if #1+ \limsup_{#2} \else {\if #1- \liminf_{#2} \else
                  \lim_{#2}\left(\begin{array}{l}\sup\\
                           \inf\end{array}\right) \fi} \fi} 
\def\cite#1{(#1)} 

\begin{document}
\title{Switched flow systems: pseudo billiard dynamics}

\author{Michael Blank\thanks{Russian Academy of Sciences,
                             and Georgia Institute of Technology;
                             e-mail: blank@iitp.ru},
        Leonid Bunimovich\thanks{School of Mathematics,
                                 Georgia Institute of Technology;
                                 e-mail: bunimovh@math.gatech.edu}}
\date{\today}
\maketitle

\n{\bf Abstract}. We study a class of dynamical systems which
generalizes and unifies some models arising in the analysis of
switched flow systems in manufacturing. General properties of
these dynamical systems, called pseudo billiards, as well as some
their perturbations are discussed.

\section{Introduction}

In a series of papers dynamics of some basic models in
manufacturing, called switched arrival and switched server
systems, has been studied \cite{Chase, Serrano and Ramadge 1993,
Horn and Ramadge 1997, Katzorke and Pikovsky 2000, Peters and
Parlitz 2003, Schurmann and Hoffman 1995}. These models describe
the local blocks of many systems in manufacturing and logistics,
where one production unit (e.g. a work station) has to load a
number of subsequent production units or one production unit has
to unload (serve) some foregoing production unit. Here we
introduce a new class of dynamical systems, called pseudo
billiards, which contains switched arrival and switched server
systems as well as some natural generalizations of these models.
Dynamical and statistical properties of pseudo billiards and
their stability against random and some deterministic
perturbations are studied. It is shown that pseudo billiards can
demonstrate chaotic, stable and neutral behavior. Sufficient
conditions for each type of behavior and examples with their
coexistence are given. The dynamics of the switched arrival
system is considered in Subsection 2.1 and the switched server
system in Subsection 2.3 respectively. Some perturbations of
these models are discussed in Subsection 2.2 and 2.3. %

\section{Mathematical background}

\subsection{Standard model}

Consider a system consisting of $N$ buffers, and one server (see
Fig.~\ref{sw-flow-system}). We refer to the content of a buffer
as ``material'', and depending on the setup the material will be
considered as a continuous object (e.g., fluid) or a discrete one
(e.g., packets). Material is removed from buffer $i$ at a fixed
rate $\rho_i>0$ ($\sum_i^N\rho_i=1$). On the other hand, the
server delivers material to a selected buffer at unit rate, and
the selection is done as follows: a threshold $\nu_i\ge0$ is
assigned to each buffer and the server switches instantaneously
to the buffer in which the level of material falls below the
assigned threshold. Otherwise the switch between the buffers
happen only after the complete work on the current buffer is
finished. If two or several buffers become empty at the same
moment of time then server chooses a buffer according to some
(cyclic) order.
To avoid overflow or complete draining of all buffers a balance
condition is necessary: the rate of filling is equal to
$\sum\rho_i$.

\Bfig(320,150)
      {\bline(0,0)(1,0)(150)   \bline(0,0)(0,1)(150)
       \bline(0,150)(1,0)(150) \bline(150,0)(0,1)(150)
       \thicklines
       \bline(55,10)(1,0)(40)  \bline(55,25)(1,0)(40)
       \bline(55,10)(0,1)(15)  \bline(95,10)(0,1)(15)
       \put(62,15){server}
       \bline(10,60)(1,0)(20)  \bline(10,60)(0,1)(40) \bline(30,60)(0,1)(40)
       \bline(40,60)(1,0)(20)  \bline(40,60)(0,1)(40) \bline(60,60)(0,1)(40)
       \put(80,70){{\Large\bf\dots}}
       \bline(120,60)(1,0)(20)  \bline(120,60)(0,1)(40) \bline(140,60)(0,1)(40)
       \bezier{50}(5,90)(35,90)(65,90)  \bezier{25}(115,90)(130,90)(145,90)
       \bezier{15}(5,65)(20,65)(35,65)  \put(17,75){1} \put(1,60){$\nu_1$}
       \bezier{15}(35,70)(50,70)(65,70) \put(47,75){2} \put(63,63){$\nu_2$}
       \bezier{15}(115,65)(130,65)(145,65) \put(127,75){N}\put(107,60){$\nu_N$}
       \thinlines
       \put(65,27){\vector(-3,2){45}}   \put(68,27){\vector(-1,2){15}}
       \put(85,27){\vector(3,2){45}}    \put(80,45){{\Large\bf\dots}}
       \put(63,40){$\sum\rho_i$}
       \put(20,93){\vector(0,1){25}}    \put(50,93){\vector(0,1){25}}
       \put(130,93){\vector(0,1){25}}
       \put(10,110){$\rho_1$} \put(40,110){$\rho_2$} \put(118,110){$\rho_N$}
       \put(65,-13){(a)}
       \put(170,0){\bpic(150,150){
       \bline(0,0)(1,0)(150)   \bline(0,0)(0,1)(150)
       \bline(0,150)(1,0)(150) \bline(150,0)(0,1)(150)
       \thicklines
       \bline(55,10)(1,0)(40)  \bline(55,25)(1,0)(40)
       \bline(55,10)(0,1)(15)  \bline(95,10)(0,1)(15)
       \put(62,15){buffer}
       \bline(10,60)(1,0)(20)  \bline(10,60)(0,1)(40) \bline(30,60)(0,1)(40)
       \bline(40,60)(1,0)(20)  \bline(40,60)(0,1)(40) \bline(60,60)(0,1)(40)
       \put(80,70){{\Large\bf\dots}}
       \bline(120,60)(1,0)(20)  \bline(120,60)(0,1)(40) \bline(140,60)(0,1)(40)
       \bezier{50}(5,90)(35,90)(65,90)  \bezier{25}(115,90)(130,90)(145,90)
       \put(17,75){1} \put(47,75){2} \put(127,75){N}
       \thinlines
       \put(20,57){\vector(3,-2){45}}   \put(53,57){\vector(1,-2){15}}
       \put(130,57){\vector(-3,-2){45}} \put(80,45){{\Large\bf\dots}}
       \put(63,40){$\sum\rho_i$}
       \put(20,118){\vector(0,-1){25}}    \put(50,118){\vector(0,-1){25}}
       \put(130,118){\vector(0,-1){25}}
       \put(10,110){$\rho_1$} \put(40,110){$\rho_2$} \put(118,110){$\rho_N$}
       \put(65,-13){(b)}
       }}
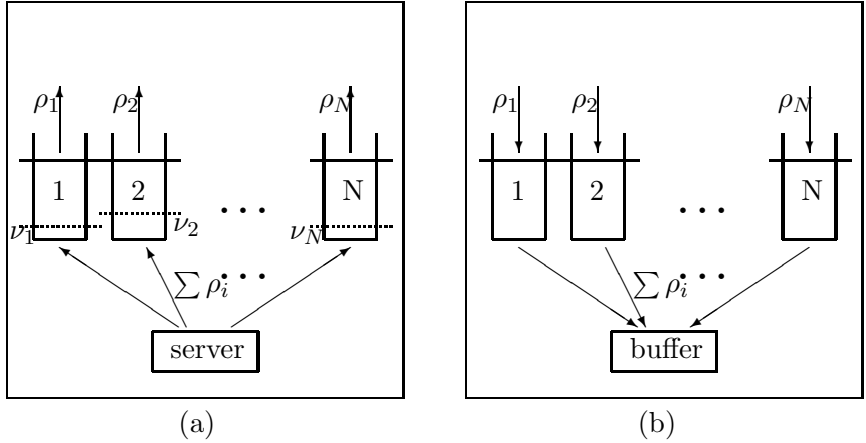
      }{(a) Switched arrival system; (b) switched server system
        \label{sw-flow-system}}

The state of this system is completely described by the
$N$-dimensional vector $x\in\IR^N$ whose $i$-th coordinate $x_i$
represents the amount of material in the $i$-th buffer and the
number of the buffer serving at the moment. Thus the evolution
can be described as follows. The current state (which we
associate to a point in $\IR^N$) moves uniformly and linearly
inside the bounded region $\{x\in\IR^N: ~ \sum_ix_i=\sum_i\rho_i,
~ \nu_i \le x_i \le \bar\nu_i\}$, i.e.
$$\dot x(t):=e^{(i)}-\rho ,$$
where $e^{(i)}$ is the $i$-th canonical unit vector in $\IR^{N}$
and $\rho:=(\rho_{1},\dots,\rho_{N})\in\IR^{N}$. After a collision
with the boundary of this region, i.e. when one of the
coordinates, say $x_{j}$, becomes equal to $\nu_{j}$, the velocity
changes instantaneously to a new vector $e^{(j)}-\rho$ which
depends only on the number of the buffer $j$ which will be served.

In the trivial cases $N=1$ or $N=2$ the dynamics is, obviously,
periodic. In the simplest nontrivial case, when $N=3$,
$\rho_i=1/3$, and $\nu_i=0$ the dynamics of this system is equivalent
to the motion of a particle in the equilateral unit triangle (see
Fig.~\ref{billiard}(a)) which moves with the constant velocity
inside the triangle and changes it instantaneously at the moment
of a collision with the boundary to the velocity perpendicular to
the corresponding side of the triangle. The motion inside the
triangle corresponds to the situation when the server delivers
material to the selected buffer, while switching between
different buffers is described by the collision rule. Due to the
symmetry the direction of the velocity after the collision is
perpendicular to the side of the triangle where the collision
took place.

Due to similarity of this dynamics to the motion of a ball in a
billiard game, we call systems of this type pseudo billiards (in
the literature one can find also the notion ``strange billiard''
used to describe a special particular case of pseudo billiards,
see, e.g. \cite{Peters and Parlitz 2003, Schurmann and Hoffman
1995}).

To be precise, by the {\em pseudo billiard} we shall mean
dynamics of the particle which moves with the constant velocity
inside a given region (not necessarily a polyhedron) and changes
it instantaneously at the moment of a collision with the boundary
to the velocity defined by a given vector field (not necessarily
a constant one) on the boundary of the region. Note that a
similar situation arises for billiards in a strong magnetic or in
the gravitational field, where only the angle with the field
matters.

A straightforward generalization of the considered above case of
three equal buffers to the general system of buffers with $N\ge3$,
$\rho_i=1/N$, and $\nu_i=0$ leads to the motion of the particle
with unit velocity in the equilateral unit $(N-1)$-dimensional
pyramid up to a moment of a collision with some
$(N-2)$-dimensional face of the pyramid when the direction of the
velocity changes instantaneously to the one perpendicular to this
face.

More general model with arbitrary rates $\rho_i>0$ and small
enough thresholds $\nu_i\ge0$ leads to a similar picture, where
the $i$-th vertex of the equilateral unit $(N-1)$-dimensional
pyramid is ``cut'' by the hyperplanes perpendicular to the
corresponding axes at the level $\nu_i$ and the ``perpendicular
reflection law'' is changed to the law defined by the relation
between $\rho_i, ~ 1\le i\le N$ in the case of the reflection law
from the `main' part of the face, while the reflection from the
``cut'' part is inherited from the corresponding `normal' face
(see details in Section~\ref{S:pert}).

Therefore all these situations can be considered as special cases
of the motion of the particle with unit velocity in a polyhedron
$\cX\in\IR^d$ with the ``reflection law'' on boundaries defined
by a given vector field. We will refer to the simplest situation
when $\cX\in\IR^d$ is the equilateral unit $d$-dimensional
pyramid and the restriction of the vector field to each of its
$(d-1)$-dimensional faces is a constant vector-field as to the
{\em standard model}.

Denote by $X$ the union of $(d-1)$-dimensional faces of $\cX$ and
by $\Sigma$ the union of its $(d-2)$-dimensional faces. Note that
the set $\Sigma$ has zero $(d-1)$-dimensional Lebesgue measure.
The {\em first return (Poincare)} map from the set $X$ into
itself we denote by $\map$. Recall that for any point $x\in X$
its image $\map x$ is defined as the next intersection of the
trajectory of the point $x$ with $X$.

\Bfig(320,150)
      {\bline(0,0)(1,0)(150)   \bline(0,0)(0,1)(150)
       \bline(0,150)(1,0)(150) \bline(150,0)(0,1)(150)
       \bline(50,50)(1,0)(80)  \bline(50,50)(0,1)(80)
       \bline(50,50)(-1,-1)(40)
       \thicklines
       \bline(15,15)(1,3)(35)  \bline(15,15)(3,1)(105)
       \bline(50,120)(1,-1)(70)
       \thinlines
       \put(65,32){\vector(0,1){15}}
       \put(32,67){\vector(1,0){15}}
       \put(85,85){\vector(-1,-1){15}}
         \bline(50,100)(1,0)(20)  \bline(50,100)(-1,-1)(10)
         \bline(70,100)(-3,-1)(30)
       \put(65,-13){(a)}
       \put(170,0){\bpic(150,150){
       \bline(0,0)(1,0)(150)   \bline(0,0)(0,1)(150)
       \bline(0,150)(1,0)(150) \bline(150,0)(0,1)(150)
       \thicklines
       \bline(10,15)(1,0)(130)  \bline(10,15)(1,3)(43.25)
       \bline(140,15)(-2,3)(86.5)
       \thinlines
       \bline(10,15)(3,1)(106.5)  \bline(140,15)(-4,5)(92)
       \bline(53,145)(1,-4)(32.5)
       \put(65,-13){(b)}
       }}
      }{(a) Pseudo billiard; (b) conditions of chaoticity \label{billiard}}

There are two major approaches to a rigorous analysis of dynamical
systems of this type. One is based on construction of the so
called Markov partition and study of the corresponding symbolic
dynamics. The second, operator approach, deals with the action of
the dynamical system in the space of measures. Both of these
approaches were already used for the analysis of the standard
model (at least of its 2-dimensional version, see, e.g.,
\cite{Chase, Serrano and Ramadge 1993, Katzorke and Pikovsky
2000}), and we will start with general discussion of their
advantages and restrictions.

Consider first a Markov partition. Observe that in the simplest
case $\nu_i=0$ the dynamical system $(\map,X)$ has the property
that each $(d-1)$-dimensional face is mapped onto the union of
some other $(d-1)$-dimensional faces and this mapping is
one-to-one. A slight generalization of this property leads to
partition of the phase space $X$ into a collection of connected
components $X_i$ with piecewise smooth boundaries satisfying the
following conditions:\par%
\begin{itemize}
\item $\cup_i X_i=X$,
\item ${\rm Int}(X_i)\cap {\rm Int}(X_j)=\emptyset$ ~ $\forall i\ne j$,
\item $\map X_i = \cup_{j\in I_i}X_j$,
\item $\map|_{X_i}$ is a homeomorphism for each $i$,
\end{itemize}
which we shall call the {\em strong Markov property}, and the
corresponding partitions will be called strong Markov partitions.
The word `strong' is used here to underline the distinction to the
usage of Markov partitions in the case of hyperbolic systems, where
only ``stable'' boundaries of elements of the partition are preserved.


Consider the action of the map $\map$ in the space ${\cal M}(X)$
of probability measures on $X$. Recall that a measure $\mu\in{\cal
M}(X)$ is called $\map$-{\em invariant} if for any measurable set
$A\subseteq X$ we have $\mu(A)=\mu(\map^{-1}A)$, and {\em
ergodic} with respect to $\map$ if $\mu$-measure of any
$\map$-invariant set is equal either to 0 or to 1. A dynamical
system $(\map, X)$ with a nonergodic invariant measure $\mu$ can
be uniquely (up to sets of $\mu$-measure zero) decomposed into
{\em ergodic components}, i.e. such subsets $\{Y_{i}\}$ that for
each $i$ the restriction of the measure $\mu$ to $Y_{i}$ is an
ergodic invariant measure for the restriction of the map $\map$
to $Y_{i}$.

Let us fix a reference measure $m$ on $X$. Since in the examples
considered in the paper the phase space typically will be a
subset of the $d$-dimensional Euclidean space $\IR^{d}$ and only
Lebesgue reference measures will be considered. We shall say that
a dynamical system $(\map, X)$ is {\em strongly chaotic} if there
is a $\map$-invariant measure $\mu$ absolutely continuous with
respect to the reference measure $m$ (i.e. $\mu$ has a density
with respect to $m$).

The simplest maps having this property are the so called
expanding maps. A map $\map$ is {\em expanding} if its Jacobian
$D\map|_{x}$ is an expanding matrix at any point $x\in X$. Recall
that a matrix $A$ is called {\em expanding} if there is a number
$\la_A>1$ called {\em expanding constant} such that for any
vector $z$ we have $|Az|\ge\la_A|z|$, where $|z|$ is the
Euclidean norm of the vector $z$.

\begin{theorem} \label{t:markov} 
Let $\map:X\to X$ be a piecewise $C^1$ map having a finite strong
Markov partition $\{X_i\}_{i=1}^K$. Then the phase space of the
dynamical system $(\map,X)$ can be decomposed into at most $K$
ergodic components and in each ergodic component all Lyapunov
exponents are either all positive or all equal zero. If $d=2$ and
the map $\map$ is continuous, then the system is either ergodic
and strongly chaotic (has an a.c.i.m.), or pure neutral (each
trajectory is eventually periodic).
\end{theorem}


\proof Elements of the Markov partition $\{X_i\}$ can be divided
into groups (components of transitivity) with respect to the
possibility to go from one of them to another under the action of
$\map$. Namely two elements $X_{i},X_{j}$ of the Markov partition
belong to the same component of transitivity if there is a
positive integer $n$ such that the interior part ${\rm
Int}(\map^{n}X_{i}\cap X_{j})\ne\emptyset$.

Consider a sequence of matrices $A_n:=\prod_{k=0}^{n-1}
D\map(\map^{k} x)$, where $D\map(y)$ is the Jacobian of the map
$\map$ at point $y$. According to Multiplicative Ergodic Theorem
(see, e.g., \cite{Katok and Hasselblatt 1995})
for a.a. $x\in X$ the limits as
$n\to\infty$ of logarithms of eigenvalues of the matrix $A_n$
divided by $n$ (called Lyapunov exponents of the map $\map$)
exist and do not depend on $x$ (on each component of
transitivity). Observe now that due to the assumption on existence
of a finite strong Markov partition all eigenvalues should be
nonnegative (otherwise we will get the contraction, which
contradicts the strong Markov property).

If some of the Lyapunov exponents are equal to 0 then the
corresponding components of transitivity can be nonergodic
(consider, for example, a rational rotation of the circle $\map
x=x+\alpha \pmod{2\pi}$ where $\alpha/\pi$ is a rational number).
On the other hand, on each ergodic component the system is either
pure expanding, or pure neutral.

Consider now the case $d=2$ and assume that the only Lyapunov
exponent of this system is positive on some component of
transitivity. Since we have a finite Markov partition, there is a
constant $n_0$ such that the map $\map^{n_0}$ is strictly
expanding. Clearly the map $\map^{n_0}$ is still continuous and
the strong Markov property holds for it with respect to the same
Markov partition. Thus due to the continuity of the map $\forall
i$ ~ $\map^{n_0}X_i$ is a union of {\em neighboring} elements of
the Markov partition, and the Lebesgue measure of this union is
strictly larger than the Lebesgue measure of the set $X_i$.
Therefore there is another constant $n_1$ such $\forall i$ ~
$\map^{n_0\cdot n_1}X_i=X$. Thus the map $\map^{n_0\cdot n_1}$ is
uniformly expanding and has an absolutely continuous (with
respect to Lebesgue measure) invariant measure. \qed

If the dimension $d>2$ the situation can be more complex and a
chaotic regime may coexist with a neutral periodic one. Moreover
if we drop the assumption of the continuity of the map $\map$ the
coexistence may take place even if $d=2$ and a map is piecewise
linear (see below). In the multidimensional case it is much easier
to construct a counterexample of the absence of ergodicity of the
system because one can consider a direct product of a chaotic
$2$-dimensional system and a neutral one, which is certainly
nonergodic. On the other hand, it turned out to be rather
difficult (but still possible) to construct a counterexample with
a map having discontinuities only on the boundary of a strong
Markov partitions and being non ergodic. The following result
demonstrates coexistence of chaotic and neutral regimes in a
2-dimensional pseudo billiard.

\begin{theorem} \label{t:coexistence-markov} 
There exists a 2-dimensional pseudo billiard generated by a
discontinuous map $\map$ such that it has both strongly chaotic
and neutral components of transitivity.
\end{theorem}

\proof Consider the pseudo billiard depicted in
Fig.~\ref{billiard-coe}. Here a chaotic regime (on the component
AG,FG,FE,BC,CD) coexists with a neutral one (on another component
AB,ED). In this construction we assume that the pairs of straight
lines (AB, GD), (BG,CE), (AE,BD), and (AE,GF) are parallel to
each other, and that $|GF|<|FE|$, $|CD|<|BC|<|GF|+|EF|$, and
$|AG|<|BC|+|CD|$, which can be easily achieved. (Observe that we
can weaken these assumptions by assuming that FE is mapped onto
the union of AG and GF and the statement about the coexistence of
two components still remains valid.) Under these assumptions the
map is piecewise expanding on the components AG,FG,FE,BC,CD, and
it is neutral on the remaining part of the boundary. The
invariance of both these components follows from the
construction. \qed

\Bfig(320,150)
      {\bline(0,0)(1,0)(200)   \bline(0,0)(0,1)(150)
       \bline(0,150)(1,0)(200) \bline(200,0)(0,1)(150)
       \thicklines
       \bline(15,15)(1,0)(117) 
       \bline(15,15)(0,1)(70)  
       \bline(15,85)(1,4)(7)   
       \bline(22,113)(1,1)(24) 
       \bline(132,15)(2,3)(32.5) 
       \bline(150,85)(2,-3)(14.5) 
       \bline(150,85)(-2,1)(104) 
       \put(75,15){\vector(1,4){5}} 
       \put(95,113){\vector(-1,-4){5}} 
       \put(15,45){\vector(1,0){15}} 
       \put(150,43){\vector(-3,2){15}} 
       \put(18,97){\vector(2,3){10}} 
       \put(34,125){\vector(-1,-4){5}} 
       \put(157,74.5){\vector(-1,-4){5}} 
       \thinlines
       \bezier{50}(15,85)(74,45)(132,15)
       \bezier{50}(46,136)(104,101)(161,65)
       \bezier{50}(15,15)(30,75)(46,136)
       \bezier{50}(150,85)(82,85)(15,85)
       \bezier{50}(150,85)(141,50)(132,15)
       \put(10,5){A}   \put(130,5){B} 
       \put(5,83){G}   \put(15,117){F}
       \put(36,137){E} \put(166,63){C}  \put(152,82){D}
       \put(220,0){\bpic(100,150){
       \bline(0,0)(1,0)(100)   \bline(0,0)(0,1)(150)
       \bline(0,150)(1,0)(100) \bline(100,0)(0,1)(150)
       \thicklines
       \put(43,130){AG}
       \put(50,125){\vector(-1,-1){15}} \put(50,125){\vector(+1,-1){15}}
       \put(23,97){BC}  \put(60,97){CD} \put(57,102){\vector(-1,0){15}}
       \put(30,92){\vector(-1,-1){15}}  \put(30,92){\vector(+1,-1){15}}
       \put(3,66){FE}   \put(40,66){GF} \put(37,70){\vector(-1,0){15}}
       \put(10,76){\vector(0,1){57}}    \put(10,133){\vector(1,0){30}}
       \put(23,20){AB}  \put(60,20){ED} \put(57,25){\vector(-1,0){15}}
                                        \put(42,25){\vector(1,0){15}}
       \thinlines
       }}
      }{Coexistence of chaotic and neutral regimes \label{billiard-coe}}

The main disadvantage of the approach based on Markov partitions
is that the construction of a partition is rather unstable with
respect to various perturbations (except perturbations adapted to
dynamics in a sense that the Markov partition is preserved under
the perturbation, which is a very rare situation) and that
`typically' a dynamical system does not have a finite Markov
partition (especially a strong one). Moreover, even if a
dynamical system has a finite Markov partition, then `typically'
under arbitrary small perturbations this property does not
survive.

Consider the simplest case of the standard model with $N=3$ and a
small `cutting' of one of its faces. Here by the cutting
operation we mean the partition of the original phase space in
two parts by a hyperplane. Choosing one of those parts we get the
phase space for the perturbed system, and by the amplitude of the
perturbation (its smallness) we mean the volume of the additional
part. It is straightforward to show that for almost all
`cuttings' this new system does not have a strong finite Markov
partition. To deal with this situation we shall apply the
operator approach.

Observe that in this case some iteration of the map $\map$
restricted to $X$ is a {\em piecewise expanding} (PE) map, which
means that there is a collection of connected sets $X_i$ with
piecewise smooth boundaries satisfying the properties:
\begin{itemize}
\item $\cup_i X_i=X$,
\item ${\rm Int}(X_i)\cap {\rm Int}(X_j)=\emptyset$ ~ $\forall i\ne j$,
\item $\map|_{X_i}$ is a $C^1$-homeomorphism for each $i$,\par%
\item the Jacobian $D(\map|_{X_i})$ is an expanding matrix for all
      $i$ and $x\in X_{i}$.
\end{itemize}

The operator approach is well developed for the case of piecewise
expanding maps (see, e.g., \cite{Blank 2001}). For instance, it
can be shown that if a map is piecewise expanding, then it is
strongly chaotic. Recall that a map $\phi: X\to X$ is called an
{\em absolutely continuous homeomorphism} if it is one to one and
for any set $A$ of Lebesgue zero measure its image $\phi(A)$ has
also zero Lebesgue measure.

The following statement is one of the basic observations in the
theory of piecewise expanding maps.

\begin{theorem} \label{t:conj}
If a continuous map $\map$ is conjugated to a piecewise expanding map
$\tilde \map$ by an absolutely continuous homeomorphism $\phi$
(i.e. $\phi\map = \tilde\map\phi$), then all Lyapunov exponents of the
map $\map$ are positive and $\map$ has an absolutely continuous
invariant measure $\mu_\map$.
\end{theorem}

Consider now the simplest standard model, i.e. the case when $\cX$
is the equilateral unit pyramid with triangular
$(d-1)$-dimensional faces and on each face the vector field is
perpendicular to it.

\begin{theorem}\label{t:standard}
If $d>2$ then for the standard model, the first return map $\map$
is not piecewise expanding. On the other hand, the $(d-2)$-th
iteration of this map (i.e. $\map^{d-2}$) is piecewise expanding,
and thus this map is strongly chaotic.
\end{theorem}

\proof Let $d=3$. Fix one of the $(d-1)$-dimensional faces of the
pyramid and denote by $A,B,C$ its vertices and by $O$ the
orthogonal projection of the opposite vertex to this face. Choose
two points $D\in AO$ and $E\in BO$ such that the lines $AB$ and
$DE$ be parallel to each other. Then the image of the interval
$DE$ under the action of $\map$ is again an interval of the same
length parallel to $DE$. Therefore there is no expansion in this
direction.

Observe that if $d=3$ then the described above neutral direction
cannot be mapped again to the neutral direction. Therefore in the
case $d=3$ the second iterate of our map becomes expanding.

If the dimension $d\ge3$, on each of the $(d-1)$-dimensional faces one
can construct a $(d-2)$ dimensional simplex such that under the action
of $\map$ it will be mapped isometrically to the $(d-2)$ dimensional
simplex on another face of $\cX$. However, similarly to the
$2$-dimensional case, only a $(d-3)$-dimensional face of this image
will be mapped isometrically under the action of $\map^2$. Thus after
each successive iteration the dimension of the isometrically mapping
simplex is decreasing by 1, which finishes the proof. \qed

If the assumption that the set of end-points of boundary segments
is mapped into itself is not satisfied, the system can demonstrate
any type of behavior, starting from the chaotic one and up to the
appearance of globally stable periodic trajectories as well as to
their coexistence.

In this case one can expect a complex behavior even for the
standard model. The following simple assertion describes
necessary and sufficient conditions for the existence of stable
and chaotic regimes for this model.

\begin{theorem}\label{t:st.model-chaos}
Consider the standard model. If for each vertex of $\cX$ the
straight line drawn parallel to the constant vector-field
corresponding to the opposite face intersects the face in its
inner point (see Fig.~\ref{billiard}(b)) then this pseudo billiard
system is strongly chaotic. Otherwise a pseudo billiard can have
a globally stable periodic trajectory.
\end{theorem}

\proof It is straightforward to check that under the assumptions
of this theorem the map $\map$ satisfies the strong Markov
property with the Markov partition consisting of pieces with
piecewise linear boundaries of $(d-1)$ dimensional faces of $\cX$.
Indeed, if the straight line from a vertex $A$ drawn parallel to
the constant vector-field corresponding to the opposite face
intersects this face, then the image of this face coincides with
the union of $(d-1)$ dimensional faces having a common vertex
$A$. Thus, the preimages of the $(d-1)$ dimensional faces of
$\cX$ generate the strong Markov partition. Observe now that the
map $\map$ is continuous, linear on elements of the partition, and
that some iterate of each element of the partition covers $\cX$.
Therefore there exists a positive integer $n_{*}$ such that the
map $\map^{n_*}$ is piecewise expanding and maps each element of
its strong Markov partition onto the complete phase space, and
thus the system is strongly chaotic.

To show that if our condition does not hold the system can have a
regular dynamics, consider the 2-dimensional case and assume in
the triangle $ABC$ the face $AB$ is mapped to $AC$, the face $AC$
is mapped to $AB$ and $BC$ is mapped to $AB\cup AC$. Observe that
it occurs if the straight lines from vertices $B$ and $C$ drawn
parallel to the constant vector-field corresponding to the
opposite faces (namely $AC$ and $AB$) do not intersect them. It
is easy to show that in this situation the point $A$ is a globally
stable fixed point of the map $\map$. \qed

\begin{remark} At the first sight it looks like that the dynamics of
this type can be analyzed (at least in the 2-dimensional case) by
making use of the idea that the map $\map$ is topologically
conjugated to the binary one (i.e. to the map $x\to 2x\pmod1$).
Observe though that the topological conjugation is not enough to
deduce any result about metric properties of the map (like the
existence of a.c.i.m., etc.). In fact, one needs to have at least
an absolutely continuous conjugation. Besides, even in the
2-dimensional case this map is topologically conjugated to the
map $x\to1-\{2x\}$ (where $\{\cdot\}$ stays for the fractional
part of a number) rather than to the binary map. It is not known
whether or not such a conjugation exists in the multidimensional
case, which rules out a generalization of this approach to $N>3$
(recall that the dimension of the phase space of the first return
map is $N-2$).
\end{remark}

\subsection{Perturbations}\label{S:pert}

It is reasonable to consider the case with nonzero thresholds
$\nu_{i}$ as a perturbation of the standard model \cite{Peters
and Parlitz 2003}. We shall consider here three types of
perturbations: ``cutting'' (describing nonzero thresholds), time
discretization (corresponding to the case when all operations are
performed in discrete time), and space discretization (when the
work is a collection of discrete objects consisting, e.g., of
finite size packets).

The ``cutting'' perturbation is defined by means of a collection
of hyperplanes in $\IR^{d}$ intersecting the polyhedron $\cX$. A
convex region obtained from $\cX$ after the intersection defines
the phase space of the perturbed system. The ``cut'' map is
constructed as follows: a point moves along the corresponding
vector-field until it reaches either the boundary of $\cX$ or the
``cut'' piece. In the first case the map was already defined,
while in the second case the vector-field (on the ``cut'' piece)
is inherited from the point of the boundary of $\cX$ to where the
trajectory would come without the ``cutting''. Namely, one should
prolong the trajectory up to the intersection with the boundary
and the vector field at the intersection point defines the new
vector-field (see Fig.~\ref{f:cutting-def}). Let $\xi$ be a point
on the boundary of $\cX$ and let $\xi', \xi''$ be its images due
to ``cutting''. Then the vector-field at the point $\xi'$ is the
same as the one at the point $\xi''$. Here by images of a point
due to cutting we mean intersections of a trajectory started from
this point with the old and new boundaries respectively. Observe
that the image of the point $\xi$ is the point $\xi''\in X$, the
line $\xi\xi'$ is parallel to $BF$ and $\xi'\xi'''$ is parallel
to $AD$, while the directions $FB$ and $DA$ correspond to the
restrictions of the vector-field to $AC$ and $BC$ respectively.


\Bfig(320,150)
      {\bline(0,0)(1,0)(150)   \bline(0,0)(0,1)(150)
       \bline(0,150)(1,0)(150) \bline(150,0)(0,1)(150)
       \thicklines
       \bline(10,15)(1,0)(130)  \bline(10,15)(1,2)(65)
       \bline(140,15)(-1,2)(65)
       \thinlines
       \bline(10,15)(5,1)(118)     
       \bline(140,15)(-5,1)(118)   
       \bline(75,145)(-1,-4)(32.5) 
       \put(5,5){A} \put(140,5){C} \put(64,140){B}
       \put(7,33){E} \put(37,5){F} \put(130,39){D} \put(122,55){G}
       \bline(60,15)(1,4)(26.5)    \put(57,5){$\xi$}
             \put(63,61){$\xi'$}   \put(90,120){$\xi''$}  
       \bline(114,66.5)(-5,-1)(87) \put(115,64){$\xi'''$}
       \thicklines
       \bline(118,58)(-1,0)(87) \put(21,55){I} 
       \put(65,-12){(a)}
      \put(170,0){\bpic(150,150){
       \bline(0,0)(1,0)(150)   \bline(0,0)(0,1)(150)
       \bline(0,150)(1,0)(150) \bline(150,0)(0,1)(150)
       \bline(15,15)(1,0)(120) \bline(15,15)(0,1)(120)
       \bline(15,135)(1,0)(120) \bline(135,15)(0,1)(120)
       \bezier{50}(55,15)(55,75)(55,135) \put(52,5){A}
       \bezier{50}(95,15)(95,75)(95,135) \put(93,5){B}
       \bezier{50}(15,55)(75,55)(135,55) \put(7,5){C}  \put(132,5){C}
       \bezier{50}(15,95)(75,95)(135,95) \put(3,50){A} \put(3,90){B}
       \bezier{75}(15,15)(75,75)(135,135)
       \bezier{50}(38,15)(38,75)(38,135) \put(35,5){F}
       \thicklines
       \bline(55,55)(-2,5)(16) \bline(15,135)(3,-5)(24)
       \bline(55,55)(1,-3)(13.3)  \put(66,5){E}
       \bezier{50}(68,15)(68,75)(68,135)
       \bline(95,95)(-2,3)(26.5) \bline(95,95)(4,-5)(32)
       \bline(135,15)(-1,5)(8)
       \bezier{50}(127,15)(127,75)(127,135) \put(123,5){D}
       \bline(55,35)(1,-1)(10)  \bline(55,71)(-1,1)(10)
       \put(65,-12){(b)}
      }}
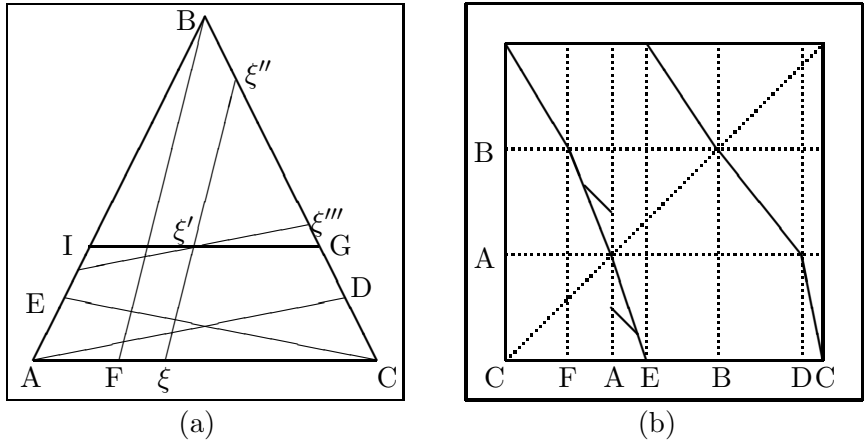
      }{(a) Cutting perturbation, (b) modifications are shown near the point A
         \label{f:cutting-def}}


\begin{theorem}\label{t:standard-cutted}
Let the standard model $(\map,\cX)$ satisfy the assumptions of
Theorem~\ref{t:st.model-chaos}. Suppose that the restriction of
the vector-field to each $(d-1)$-dimensional face is orthogonal
to it. Consider such ``cutting'' that under the action of $\map$
each of the new $(d-1)$-dimensional faces obtained as a result of
``cutting'' is mapped to the same $(d-1)$-dimensional face as
their ancestral faces do. Then the resulting system is strongly
chaotic.
\end{theorem}

\proof Observe that the ``cut'' system is again piecewise
expanding and is piecewise linear on domains consisting of a
finite number of simplexes. It is known \cite{Tsujii 2001} that
under these assumptions the map has an absolutely continuous
(with respect to Lebesgue measure) invariant measure, i.e. it is
strongly chaotic. \qed

Theorem~\ref{t:standard-cutted} gives only sufficient conditions
for the strong chaoticity. However they are close to the necessary
ones. In particular, one cannot replace the assumption that the
restriction of the vector-field to each $(d-1)$-dimensional face
is orthogonal to it by the assumption that the map $\map$ is
piecewise expanding. The following example (see
Fig.~\ref{f:contr}) shows that this cannot be done even if $d=2$.

\begin{lemma} The map corresponding to the ``cut'' two-dimensional
system depicted in Fig.~\ref{f:contr} has a stable periodic
trajectory.
\end{lemma}

\Bfig(150,150)
      {\bline(0,0)(1,0)(150)   \bline(0,0)(0,1)(150)
       \bline(0,150)(1,0)(150) \bline(150,0)(0,1)(150)
       \thicklines
       \bline(10,15)(1,0)(130)  \bline(10,15)(1,2)(65)
       \bline(140,15)(-1,2)(65)
       \thinlines
       \bline(10,15)(5,1)(118)     
       \bline(140,15)(-5,1)(118)   
       \bline(75,145)(-1,-4)(32.5) 
       \put(5,5){A} \put(140,5){C} \put(64,140){B}
       \put(7,33){E} \put(37,5){F} \put(130,39){D}
       \bline(42,15)(5,1)(89) \put(135,30){J} 
       \bline(22,38.5)(5,1)(96) \put(122,55){G} 
       \thicklines
       \bline(118,58)(-1,0)(87) \put(21,55){I} 
       \bline(131,33)(-1,-2)(9) \put(122,5){K} 
       \put(60,60){L}   \put(115,19){M}
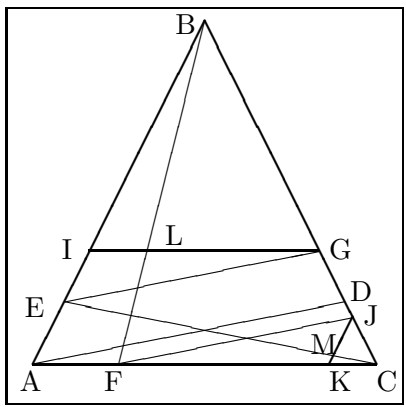
      }{Existence of a stable periodic trajectory under ``cutting'' \label{f:contr}}

\proof It is easy to see that the system before ``cuttings'' was
piecewise expanding. The lines $IG$ and $JK$ describe the
``cuttings'', which are constructed so that the conditions of
Theorem~\ref{t:standard-cutted} are satisfied. Denote the
intersection points of straight lines $IG$ with $BF$ and $JK$
with $CF$ by $L$ and $M$ respectively. Then the segment $LG$ gets
mapped under the action of $\map$ to the segment $JM$, while the
segment $JM$ gets mapped backward to the segment $LG$. Moreover,
it is not hard to show that the absolute value of the derivative
of the map along those two segments is strictly smaller than 1,
which yields the contraction. \qed

\begin{remark}
One can continue the construction in the proof above to get
arbitrary small ``cutting'' which still leads to the existence of
stable periodic trajectories. To do this let us draw a line from
the point $I$ parallel to $AD$, whose intersection with $BC$
gives a new position of the ``cutting'' with the same properties.
\end{remark}

The following Lemma gives geometric conditions equivalent to the
assumptions of Theorem~\ref{t:standard-cutted}.

\begin{lemma}\label{l:cutting}
The following two conditions are equivalent:
\begin{itemize}
\item under the action of $\map$ each $(d-1)$-dimensional face obtained
   by ``cutting'' gets mapped to the same $(d-1)$-dimensional face as
   its ancestral face;
\item the ``cut'' piece intersects with only one of the straight
   lines, drawn from vertices parallel to the constant vector-fields
   corresponding to the opposite faces.
\end{itemize}
\end{lemma}

\proof The proof is straightforward. (See
Fig.~\ref{f:cutting-def}.(a).) \qed

In view of the billiard interpretation a time discretization of
the dynamics means that inside the polyhedron $\cX$ the particle
moves in discrete steps with a given constant velocity. The only
difference occurs near the boundary where in order not to go
through the boundary (which is not permitted) the modulus of the
velocity should be changed. Thus the Poincare section of the time
discretized system at the boundary does not differ from the
original one and therefore the invariant measure does not change
under the time discretization.

\smallskip

The situation is different in case of the space discretization.
For simplicity we shall consider only uniform space
discretizations when all packets (describing the discretization)
have the same size. In this case (again as in the previous
situation) inside the polyhedron $\cX$ the particle moves in
discrete steps with the given constant velocity, whose modulus is
equal to a size of the packet. This motion continues up to the
moment when the distance to the boundary in some direction
becomes smaller then the packets' size. Since this size cannot be
made smaller, this moment has to be considered as a ``collision''
with the boundary and the velocity of the particle has to be
changed to the one defined by the vector-field on this component
of the boundary. Thus typically a trajectory of the particle
never visits the boundary. Therefore the approach based on the
Poincare section on the boundary does not make much sense.
Numerical experiments (see, e.g., \cite{Katzorke and Pikovsky
2000}) demonstrate very sensitive dependence of the statistical
properties of the discretized system on the size of packets. To
illustrate this consider the billiard type system in the square
with the sides of length $\ell\in\IZ_{+}$ and with the unit
vector-field on each side perpendicular to it. Clearly each
trajectory of the unperturbed system is periodic with the time
period $4\ell$. Take now the uniform space discretization with
the packet size $\ep>0$. If $\ep$ is a rational number, i.e.
$\ep=p/q$ where $p,q$ are positive integers then each trajectory
of the discretized system starting from the boundary is periodic
with the period $2\ell q$. On the other hand, if $\ep$ is
irrational then each trajectory densely fills a certain region.
Hence the system has an absolutely continuous invariant measure.

\subsection{Switched server system}

It makes sense to consider the situation opposite to the standard
model (see Fig.~\ref{sw-flow-system}(b)) when there is only one
buffer and several servers (see, e.g., \cite{Horn and Ramadge 1997}
for technical details). A formal definition is as follows. Let
$\cX\in\IR^{d}$ be a convex domain. Assume that $N$ vector fields
are defined on the boundary $\partial\cX$. One can think here about
a polyhedron $\cX$ with $d-1$-dimensional faces on each of which we
fix $N$ different vector fields. We consider the dynamics of a particle
moving with unit velocity inside $\cX$ and changing the velocity
at the moments of collision with the boundary $\partial\cX$ to
one of the vectors defined at the point of collision. A choice of
one of $N$ possible vector fields is governed by a certain rule
either deterministically (e.g., we enumerate the vector fields
and each time after collision increase the number of the chosen
field by one modulo $N$), or stochastically (e.g., we choose a
vector field at random with respect to a certain probabilistic
distribution $p_{1},\dots,p_{N}$). From the dynamical systems
point of view this system can be represented as a skew product. %
  Recall that the skew product of maps $F:X\to X$ and $G:X\times Y\to
  Y$ is defined as a new map $\map:X\times Y\to X\times Y$ according
  to the formula $\map(x,y):=(F(x), G(x,y))$. %
The first map $F$ is defined on the set of all integers and says
which buffer to serve next time,
while the second one $G$ is defined by the inverse branches of the
map corresponding to the pseudo billiard system described in the
previous section.

Indeed, in the simplest nontrivial case of three servers a
Poincare section of the dynamics can be defined on the boundary
of the equilateral triangle, on each side of which we consider two
vector fields perpendicular to the remaining two sides of the
triangle. Each trajectory of the system in this case coincides
with one of the inverse branches of the dynamics of the model.
Recall that the Poincare section of the standard model is
described by a piecewise monotonic map with one-to-one monotonic
branches.

Let $\cX$ be the equilateral unit pyramid with $(d-1)$-dimensional
faces on each of which we fix $N=d$ different constant vector
fields. We assume similarly to Theorem~\ref{t:st.model-chaos}
that for each vertex of $\cX$ the straight lines drawn at this
vertex parallel to each of vector fields, corresponding to the
opposite $(d-1)$-dimensional faces, intersect with them at inner
points. This assumption, which we shall call the {\em standing
assumption} holds, e.g., if the vector fields are orthogonal to
opposite faces.

\begin{theorem}\label{th-inv-det} Consider an arbitrary deterministic
  setting. Assume that the standing assumption holds. Then any
  trajectory of this system converges exponentially fast to some
  stable periodic orbit. A number of such stable periodic orbits does
  not exceed $(d+1)^{d}$.
\end{theorem}

\proof Observe that the standing assumption implies that for a
given $(d-1)$-dimensional face and any possible choice of the
constant vector field on this face the image of this face is a
proper subset of one of the opposite faces. Therefore there exits
a metrics on the boundary $\cX$ for which the dynamics is
uniformly contractive. The latter yields the needed statement.
\qed

Consider now a probabilistic version of this system, where a
choice of the vector fields at a point $x\in\partial\cX$ is
governed by some probability distribution
$p_{1}(x),\dots,p_{N}(x)$, such that $\sum_{i=1}^Np_i(x)=1$ for
each $x$. We assume that these distributions are strictly
positive, i.e. $p_i(x)>q>0, ~ i=1,\dots,N$, and piecewise
constant, i.e. there is a partition of the phase space such
that $p_i(x)$=const for all $i=1,\dots,N$ and all $x$ belonging
to the same element of the partition. This stochastic system is
equivalent to the Markov chain defined as follows. Consider a
collection of maps $\{\map_i\}$ (which we call deterministic
components) and collection of distributions
$p_{1}(x),\dots,p_{N}(x)$. Then the probability to go in one time
step from a point $x\in X$ to a set $A\subset X$ is equal to
$\sum_i p_i(x)\cdot 1_{\map_ix\cap A}$, where $1_B$ is the
characteristic function of the set $B$. Systems of this type are
called {\em iterated function system with probabilities} (IFS). %

\begin{theorem}\label{th-inv-ran}
  Suppose that the standing assumption holds. Then the Markov
  chain defined above has the unique invariant measure. \end{theorem}

\proof The Markov chain under study is ergodic due to positivity
of transition probabilities. On the other hand, each deterministic
component $\map_i$ of this system is contractive (see the proof of
Theorem~\ref{th-inv-det}), which yields the statement of Theorem
using standard properties of contractive IFS (see, e.g.
\cite{Blank 2001[Theorem 4.2.1]} for details). \qed

%

\section*{Acknowledgements}

M.B. thanks Georgia Tech for a kind hospitality and his work was
partially supported by the CRDF and RFBR grants. L.B. was
partially supported by the NSF grant DMS - 0140165 and by the
Humboldt Foundation.

\end{document}